\documentclass[11pt]{article}
\usepackage[margin=1in]{geometry}
\usepackage{amsmath,amssymb,amsthm,mathtools}
\usepackage[hidelinks]{hyperref}
\newtheorem{theorem}{Theorem}
\newtheorem{lemma}{Lemma}[section]

\newtheorem{corollary}{Corollary}
\theoremstyle{remark}
\newtheorem{remark}[lemma]{Remark}
\newcommand{\PP}{\mathbb{P}}
\newcommand{\EE}{\mathbb{E}}
\newcommand{\T}{\mathbb{T}_n}
\newcommand{\Z}{\mathbb{Z}}
\newcommand{\one}{\mathbf{1}}
\newcommand{\trel}{t_{\mathrm{rel}}}
\DeclareMathOperator{\dist}{dist}

\title{Annihilation versus anchoring of a topologically charged pair\\ on a perforated two-dimensional surface}
\author{Patrick Romanescu}
\date{}

\begin{document}
\maketitle

\begin{abstract}
Two independent random walks of opposite topological charge move on the discrete torus
$(\Z/n\Z)^2$, which contains $k$ well-separated holes. The walks annihilate when they come
within a core distance of each other, and each walk is anchored when it first touches a hole.
We prove that the probability of annihilation before anchoring converges to
$w_0/(w_0+\sum_j w_j)$, where the weights are explicit and depend only logarithmically on the
sizes of the core and of the holes; for point-like core and holes the limit is $1/(1+k)$.
We then prove a fate theorem: the pair either terminates or ends as two opposite charges anchored on two
\emph{different} holes, with probability $(k-1)/(k+1)$ in the point-like case, and the hosting pair
$\{j,l\}$ is chosen with probability proportional to $w_jw_l$.
The proof combines an abstract competition lemma for rare sets in reversible Markov chains with
sharp two-dimensional Green's function estimates.
\end{abstract}

\section{Model and results}

\subsection{Model}
Let $\T=(\Z/n\Z)^2$, $N=n^2$, with $|\cdot|$ the Euclidean norm of the representative in
$(-n/2,n/2]^2$, and $\pi$ the uniform distribution on $\T$.
A \emph{walk} is a continuous-time simple random walk on $\T$ with total jump rate $1$; its
transition kernel is $p_u(x,y)=p_u(0,y-x)$.

Fix $k\ge 1$ and exponents $\alpha_0,\alpha_1,\dots,\alpha_k\in[0,1)$; put
$\alpha_*=\max_i\alpha_i$.
\begin{itemize}
\item \emph{Holes.} $A_j=B(c_j,r_j)=\{x:|x-c_j|\le r_j\}$ with $r_j=\lfloor n^{\alpha_j}\rfloor$ if
$\alpha_j>0$ and $r_j=0$ (a single site) if $\alpha_j=0$, $j=1,\dots,k$, and
$\dist(A_i,A_j)\ge n/\log n$ for $i\ne j$.
\item \emph{Defects.} $X$ (charge $+1$) and $Y$ (charge $-1$) are independent walks started at
$x_0,y_0$.
\item \emph{Annihilation.} With $\rho=\lfloor n^{\alpha_0}\rfloor$ if $\alpha_0>0$ and $\rho=0$ if
$\alpha_0=0$, let $T_0=\inf\{u: |X_u-Y_u|\le\rho\}$.
\item \emph{Anchoring.} $T^X_j=\inf\{u:X_u\in A_j\}$ and $T^Y_j=\inf\{u:Y_u\in A_j\}$.
\end{itemize}
The start $(x_0,y_0)$ is \emph{admissible} if $|x_0-y_0|\ge \rho+n/\log n$ and
$\dist(x_0,A_j)\wedge\dist(y_0,A_j)\ge n/\log n$ for all $j$.
Let $T=T_0\wedge\min_j(T^X_j\wedge T^Y_j)$ be the time of the first event.
Almost surely exactly one walk jumps at time $T$. If the state reached at time $T$ realises several events
simultaneously, we give priority to annihilation, then to $X$, then to the smaller index $j$.
All statements below are unaffected by this convention.

For $r\ge 0$ write $h(r)=1/\EE_\pi T_{B(0,r)}$ for a single walk started from $\pi$; by
translation invariance $\EE_\pi T_{B(c,r)}$ does not depend on $c$. Put $h_j=h(r_j)$ and
$H=\sum_{j=1}^kh_j$.

\subsection{Results}

\begin{theorem}[Competition law, rate form]\label{thm:rates}
Uniformly over admissible starts and hole positions, as $n\to\infty$,
\begin{align*}
\PP(\text{annihilation is the first event}) &= \frac{h(\rho)}{h(\rho)+H}+o(1),\\
\PP(\text{the first event is the anchoring of $X$ at $A_j$}) &= \frac{h_j}{2\,(h(\rho)+H)}+o(1),
\end{align*}
and the same holds with $Y$ in place of $X$.
\end{theorem}

\begin{theorem}[Hitting rates]\label{thm:asym}
Uniformly in $0\le r\le n^{\alpha_*}$,
\[
\EE_\pi T_{B(0,r)}=\frac{2}{\pi}\,n^2\log\frac{n}{r\vee 1}\,\bigl(1+o(1)\bigr).
\]
\end{theorem}

Set $w_i=1/(1-\alpha_i)$, $W=\sum_{j=1}^kw_j$, $Q=\sum_{j=1}^kw_j^2$.
Theorem~\ref{thm:asym} gives $h_j/h(\rho)\to w_j/w_0$, hence:

\begin{corollary}[Competition law]\label{cor:A}
$\PP(\text{annihilation first})\to w_0/(w_0+W)$. If core and holes are single sites, the limit is
$1/(1+k)$; in that case $h(\rho)=h_j$ for every $n$ by symmetry, so
Theorem~\ref{thm:asym} is not needed.
\end{corollary}

\paragraph{Fate.} After the first event the dynamics continue as follows.
If the first event is annihilation, the outcome is \emph{termination}.
If a defect is anchored at $A_i$, it stops, and the other defect continues as a walk until it first
touches $\bigcup_lA_l$, say at $A_l$. If $l=i$ the charges cancel and the outcome is
\emph{termination}; if $l\ne i$ the outcome is the \emph{sustained pair} $S_{\{i,l\}}$:
charges $+1$ and $-1$ anchored on two different holes.

\begin{theorem}[Fate of a pair]\label{thm:fate}
Uniformly over admissible starts and hole positions,
\begin{align*}
\PP(\text{termination})&\to \frac{w_0}{w_0+W}+\frac{Q}{W(w_0+W)},\\
\PP(S_{\{i,l\}})&\to \frac{2\,w_iw_l}{W(w_0+W)}\qquad(i<l),\\
\PP(\text{sustained pair})&\to\frac{W^2-Q}{W(w_0+W)}.
\end{align*}
In particular, with a single hole the pair terminates with probability $\to1$, and for $k$ point-like
holes and a point-like core $\PP(\text{sustained pair})\to(k-1)/(k+1)$.
\end{theorem}

Throughout, $C,c$ denote positive constants depending only on $k$ and $\alpha_*$, whose value may
change from line to line, and $n$ is assumed large enough. We fix
\[
g=g_n=(\log n)^{1/2},\qquad t=g\,n^2,\qquad s=g^{1/2}n^2 .
\]

\section{An abstract competition lemma}

Let $Z$ be an irreducible, reversible, continuous-time Markov chain on a finite set $S$ with
stationary distribution $\pi$, relaxation time $\trel$, and
$d(u)=\max_z\|\PP_z(Z_u\in\cdot)-\pi\|_{\mathrm{TV}}$, which is non-increasing in $u$.
Let $B_0,\dots,B_m\subset S$ be nonempty (not necessarily disjoint), $T_i$ the hitting time of
$B_i$, $T=\min_iT_i$, and $I=\min\{i:Z_T\in B_i\}$.
Let $f_i(z)=\PP_z(I=i)$ and $\bar p_i=\PP_\pi(I=i)$.
For fixed $0<s<t$ let
\[
a_i=\PP_\pi(T_i\le t),\qquad a=\sum_ia_i,\qquad
\delta=\sum_{i<l}\PP_\pi(T_i\le t,\;T_l\le t).
\]

\begin{lemma}\label{lem:comp}
If $\delta<a$, then for every $i$
\[
\Bigl|\bar p_i-\frac{a_i}{a}\Bigr|\le d(s)+\frac{\PP_\pi(T\le s)+\delta}{a-\delta}.
\]
\end{lemma}
\begin{proof}
On $\{T>t\}$ the Markov property gives $\PP_\pi(I=i\mid\mathcal F_t)=f_i(Z_t)$, and
$\EE_\pi f_i(Z_t)=\bar p_i$ by stationarity. Hence
\begin{equation}\label{eq:id}
\PP_\pi(T\le t,\,I=i)=\bar p_i-\EE_\pi[\one_{\{T>t\}}f_i(Z_t)]=\EE_\pi[\one_{\{T\le t\}}f_i(Z_t)].
\end{equation}
On $\{T\le t-s\}$ the strong Markov property at $T$ gives
$\EE_\pi[f_i(Z_t)\mid\mathcal F_T]=(P_{t-T}f_i)(Z_T)$, and since $0\le f_i\le1$ and $t-T\ge s$,
$|(P_{t-T}f_i)(Z_T)-\bar p_i|\le d(s)$. On the other hand
$\PP_\pi(t-s<T\le t)\le\PP_\pi(Z\text{ visits }\textstyle\bigcup_iB_i\text{ in }[t-s,t])=\PP_\pi(T\le s)$
by stationarity. Splitting \eqref{eq:id} according to $\{T\le t-s\}$ and $\{t-s<T\le t\}$,
\begin{equation}\label{eq:c1}
\bigl|\PP_\pi(T\le t,\,I=i)-\bar p_i\,\PP_\pi(T\le t)\bigr|\le d(s)\,\PP_\pi(T\le t)+\PP_\pi(T\le s).
\end{equation}
If $I=i$ then $T_i=T$; if $T_i\le t$ and $I\ne i$ then $T_l\le T_i\le t$ for some $l\ne i$. So
\[
a_i-\delta\le\PP_\pi(T\le t,\,I=i)\le a_i,\qquad a-\delta\le\PP_\pi(T\le t)\le a,
\]
the last by the Bonferroni inequality. Therefore
$\bigl|\PP_\pi(T\le t,I=i)/\PP_\pi(T\le t)-a_i/a\bigr|\le\delta/(a-\delta)$, and dividing
\eqref{eq:c1} by $\PP_\pi(T\le t)\ge a-\delta$ gives the claim.
\end{proof}

\begin{lemma}\label{lem:start}
For every $z\in S$ and every $i$, $|f_i(z)-\bar p_i|\le d(s)+\PP_z(T\le s)$.
\end{lemma}
\begin{proof}
$f_i(z)=\PP_z(T\le s,I=i)+\EE_z[\one_{\{T>s\}}f_i(Z_s)]
=\EE_zf_i(Z_s)+\EE_z[\one_{\{T\le s\}}(\one_{\{I=i\}}-f_i(Z_s))]$, and
$|\EE_zf_i(Z_s)-\bar p_i|\le d(s)$.
\end{proof}

\begin{lemma}[Aldous--Brown]\label{lem:AB}
For nonempty $A\subset S$ and all $u\ge0$,
$\bigl|\PP_\pi(T_A>u)-e^{-u/\EE_\pi T_A}\bigr|\le\trel/\EE_\pi T_A$. Consequently, with
$h_A=1/\EE_\pi T_A$,
\[
\PP_\pi(T_A\le u)\le(u+\trel)h_A,\qquad
\bigl|\PP_\pi(T_A\le u)-uh_A\bigr|\le\trel h_A+\tfrac12(uh_A)^2 .
\]
\end{lemma}
\begin{proof}
The first statement is \cite[Theorem~1]{AB}; see also \cite[Ch.~3, Sec.~5.4]{AF}. The rest follows
from $x-\tfrac12x^2\le1-e^{-x}\le x$.
\end{proof}

\section{Heat-kernel facts on the torus}

Let $\tilde p_u$ be the kernel of the rate-$1$ continuous-time simple random walk $W$ on $\Z^2$
started at $0$, so that $p_u(0,z)=\sum_{m\in\Z^2}\tilde p_u(z+nm)$. We use two standard inputs.
\begin{itemize}
\item[(L1)] \emph{Local limit theorem} \cite[Ch.~2]{LL}: for $u\ge1$ and all $x\in\Z^2$,
$\bigl|\tilde p_u(x)-(\pi u)^{-1}e^{-|x|^2/u}\bigr|\le Cu^{-2}$.
(The covariance of $W_u$ is $(u/2)I$.)
\item[(L2)] The coordinates of $W$ are independent rate-$\tfrac12$ walks on $\Z$ with kernel
$\tilde q_u(k)=e^{-u/2}I_{|k|}(u/2)$, which is symmetric, non-increasing in $|k|$, and satisfies
$\tilde q_u(0)\le C(1+u)^{-1/2}$.
\end{itemize}

\begin{lemma}\label{lem:hk}
There are absolute constants $C,c$ such that for all large $n$:
\begin{itemize}
\item[(H0)] $p_u(x,y)\le C/(1+u)$ for $0\le u\le n^2$;
\item[(H1)] $p_u(x,y)\le C/n^2$ for $u\ge n^2$, and also for all $u\le n^2$ when $|x-y|\ge n/8$;
\item[(H2)] $p_{n^2}(x,y)\ge c/n^2$ for all $x,y$;
\item[(H3)] $G(z):=\int_0^{n^2}p_u(0,z)\,du$ satisfies
$G(z)=\frac2\pi\log\frac{n}{|z|\vee1}+O(1)$ for $|z|\le n/4$, and $G(z)\le C$ for $|z|>n/4$.
\end{itemize}
\end{lemma}
\begin{proof}
(H0). $p_u(0,z)=\prod_{i=1,2}\sum_{m\in\Z}\tilde q_u(z_i+nm)$. Since $\tilde q_u$ is
non-increasing in $|k|$, for $|z_i|\le n/2$ and $m\ne0$ the term $\tilde q_u(z_i+nm)$ is at most
$(2/n)$ times the sum of $\tilde q_u$ over a block of $\lceil n/2\rceil$ consecutive sites, and
these blocks are disjoint for distinct $m$ of the same sign. Hence
$\sum_m\tilde q_u(z_i+nm)\le\tilde q_u(0)+4/n\le C(1+u)^{-1/2}$ for $u\le n^2$.

(H1), off-diagonal part. The moment generating function of a coordinate is
$\EE e^{\theta W^{(1)}_v}=\exp(\tfrac v2(\cosh\theta-1))\le\exp(\tfrac v2\theta^2)$ for
$|\theta|\le1$, so Chernoff's bound with $\theta=\min(1,b/v)$ gives
$\PP(W^{(1)}_v\ge b)\le\exp(-\min(b^2/(2v),\,b/2))$. The torus walk is the projection of $W$, so
$\PP_x(|X_v-x|\ge n/16)\le4\exp(-\min(b^2/(2v),b/2))$ with $b=n/32$. Let $|x-z|\ge n/8$ and
$u\le n^2$. By Chapman--Kolmogorov at time $u/2$, and since every $y$ has $|x-y|\ge n/16$ or
$|y-z|\ge n/16$, (H0) and symmetry give
\[
p_u(x,z)\le2\sup_{y,y'}p_{u/2}(y,y')\;\PP_x(|X_{u/2}-x|\ge n/16)
\le\frac{C}{1+u}\exp\Bigl(-\min\bigl(\tfrac{b^2}{u},\tfrac b2\bigr)\Bigr)\le\frac{C}{n^2}.
\]

(H1), $u\ge n^2$, and (H2). Let $\tilde z$ be the representative of $z$ with $|\tilde z|\le n/\sqrt2$.
By (L1), $\tilde p_{n^2}(\tilde z)\ge(\pi n^2)^{-1}e^{-1/2}-Cn^{-4}\ge c/n^2$, which gives (H2)
because $p_{n^2}(0,z)\ge\tilde p_{n^2}(\tilde z)$. By (H0), $\sup p_{n^2}\le C/n^2$, and for
$u>n^2$ Chapman--Kolmogorov gives $p_u(x,y)\le\sup_{w}p_{n^2}(w,y)\le C/n^2$.

(H3). Let $|z|\le n/4$ and $e_u(z)=p_u(0,z)-\tilde p_u(z)\ge0$. On the event that the torus walk
is at $z$ while $W_u\ne z$, the time $\sigma=\inf\{v:|W_v|_\infty\ge n/2\}$ is at most $u$ and
$|X_\sigma-z|\ge n/4-1\ge n/8$. By the strong Markov property and (H1),
$e_u(z)\le\sup_{v\le n^2,\,|x-z|\ge n/8}p_v(x,z)\le C/n^2$ for $u\le n^2$, so
$0\le\int_0^{n^2}e_u(z)\,du\le C$. By (L1),
\[
\int_0^{n^2}\tilde p_u(z)\,du=\int_1^{n^2}\frac{e^{-|z|^2/u}}{\pi u}\,du+O(1),
\]
the error being at most $1+\int_1^\infty Cu^{-2}du$. For $z=0$ the integral is $\frac2\pi\log n$.
For $|z|\ge1$ substitute $v=|z|^2/u$:
\[
\int_1^{n^2}\frac{e^{-|z|^2/u}}{\pi u}du=\frac1\pi\int_{|z|^2/n^2}^{|z|^2}\frac{e^{-v}}{v}dv
=\frac1\pi\log\frac{n^2}{|z|^2}+O(1),
\]
because $\int_a^1(1-e^{-v})v^{-1}dv\le1$ and $\int_1^\infty e^{-v}v^{-1}dv\le1$. This proves the
first part. If $|z|>n/4$, (H1) gives $G(z)\le n^2\cdot C/n^2$.
\end{proof}

\begin{lemma}[Mixing]\label{lem:mix}
For the single walk, $d(mn^2)\le(1-c)^m$ and $\trel\le n^2/4$. For the pair $(X,Y)$ on
$\T\times\T$, $d(mn^2)\le(1-c^2)^m$ and $\trel\le n^2/4$.
\end{lemma}
\begin{proof}
By (H2), $p_{n^2}(x,\cdot)\ge c\,\pi(\cdot)$; Doeblin's argument gives the bound on $d$, and the
pair satisfies the same minorisation with $c^2$. The generator of the walk has eigenvalues
$1-\tfrac12(\cos\frac{2\pi a}{n}+\cos\frac{2\pi b}{n})$; the smallest non-zero one is
$\tfrac12(1-\cos\frac{2\pi}n)\ge4/n^2$ because $1-\cos x\ge2x^2/\pi^2$ on $[-\pi,\pi]$. The
spectral gap of a product of independent chains is the minimum of the gaps.
\end{proof}

\section{Hitting estimates for a single walk}

Let $A=B(a_0,r)$ with $0\le r\le n^{\alpha_*}$, and $L=\log\frac{n}{r\vee1}\ge(1-\alpha_*)\log n$.
For $v\ge0$ and $x\in\T$ let $\Gamma_v(x)=\EE_x\int_0^v\one_A(X_u)\,du=\sum_{y\in A}\int_0^vp_u(x,y)du$.

\begin{lemma}[Occupation]\label{lem:occ}
For $a\in A$ and $n^2\le v\le 2t$,
$\Gamma_v(a)=|A|\bigl[\tfrac2\pi L+O(1+g)\bigr]$, and for all $v\le2t$,
$\Gamma_v(a)\le|A|\bigl[\tfrac2\pi L+C(1+g)\bigr]$.
If $\dist(x,A)\ge n/\log n$ and $v\le4t$, then $\Gamma_v(x)\le|A|\,C(\log\log n+g)$.
\end{lemma}
\begin{proof}
For $a,y\in A$ we have $|a-y|\le2r\le n/4$, so by (H3)
\[
\sum_{y\in A}G(y-a)=|A|\,\tfrac2\pi L+\tfrac2\pi\sum_{y\in A}\log\frac{r\vee1}{|a-y|\vee1}+O(|A|).
\]
Each term of the last sum is at least $-\log2$. Grouping $y$ by $|a-y|\in(r2^{-j-1},r2^{-j}]$,
$j\ge0$, the sum is at most $Cr^2\sum_j(j+1)4^{-j}\le C|A|$ when $r\ge1$ (and it vanishes when
$r=0$). By (H1), $0\le\sum_{y\in A}\int_{n^2}^vp_u(a,y)du\le C|A|(v-n^2)/n^2\le C|A|g$. This proves
the first two claims. If $\dist(x,A)\ge n/\log n$, then by (H3) $G(y-x)\le\frac2\pi\log\log n+C$ for
every $y\in A$, and the contribution of $u\in[n^2,v]$ is at most $C|A|\cdot4g$.
\end{proof}

\begin{lemma}[Hitting bounds]\label{lem:hit}
\begin{itemize}
\item[(a)] For every $x$ and $v>0$,
$\PP_x(T_A\le v)\le\Gamma_{2v}(x)/\inf_{a\in A}\Gamma_v(a)$.
\item[(b)] $\PP_\pi(T_A\le n^2)\le C/L$ and $\EE_\pi T_A\ge c\,n^2L$.
\item[(c)] If $\dist(x,A)\ge n/\log n$ and $v\le2t$, then
$\PP_x(T_A\le v)\le\kappa_n:=C(\log\log n+g)/\log n\to0$.
\end{itemize}
\end{lemma}
\begin{proof}
(a) $X_{T_A}\in A$, so by the strong Markov property
$\Gamma_{2v}(x)\ge\EE_x[\one_{\{T_A\le v\}}\Gamma_v(X_{T_A})]$.
(b) Averaging (a) over $x\sim\pi$ with $v=n^2$ and using $\EE_\pi\Gamma_{2n^2}=2|A|$ and
Lemma~\ref{lem:occ}, $\PP_\pi(T_A\le n^2)\le2/(\frac2\pi L-C)\le C_1/L$. By stationarity and the
union bound, $\PP_\pi(T_A\le mn^2)\le mC_1/L$; with $m=\lfloor L/(2C_1)\rfloor$ this is at most
$\tfrac12$, so $\EE_\pi T_A\ge\tfrac12mn^2$.
(c) Since $\PP_x(T_A\le v)$ is non-decreasing in $v$ we may take $n^2\le v\le2t$, and (a) with
Lemma~\ref{lem:occ} gives the bound $C(\log\log n+g)/(\frac2\pi L-C(1+g))$.
\end{proof}

\begin{proof}[Proof of Theorem~\ref{thm:asym}]
Write $E=\EE_\pi T_A$. By stationarity and the strong Markov property,
\[
\frac{t|A|}{n^2}=\EE_\pi\int_0^t\one_A(X_u)du=\EE_\pi\bigl[\one_{\{T_A\le t\}}\Gamma_{t-T_A}(X_{T_A})\bigr].
\]
By Lemma~\ref{lem:occ} the right-hand side is at most
$\PP_\pi(T_A\le t)\,|A|[\frac2\pi L+C(1+g)]$ and at least
$\PP_\pi(T_A\le t-n^2)\,|A|[\frac2\pi L-C(1+g)]$. Since
$\PP_\pi(t-n^2<T_A\le t)\le\PP_\pi(T_A\le n^2)\le C/L$ by stationarity and
Lemma~\ref{lem:hit}(b),
\[
\frac{g}{\frac2\pi L+C(1+g)}\le\PP_\pi(T_A\le t)\le\frac{g}{\frac2\pi L-C(1+g)}+\frac CL .
\]
As $g\to\infty$ and $g/L\to0$, this gives $\PP_\pi(T_A\le t)=\frac{\pi g}{2L}(1+o(1))$, uniformly
in $r$. By Lemma~\ref{lem:hit}(b), $t/E\le g/(cL)\to0$, and by Lemmas~\ref{lem:AB} and
\ref{lem:mix}, $|\PP_\pi(T_A\le t)-t/E|\le\frac{t}{E}\bigl(\frac1{4g}+\frac{t}{2E}\bigr)$.
Hence $t/E=\frac{\pi g}{2L}(1+o(1))$, i.e.\ $E=\frac2\pi n^2L(1+o(1))$.
\end{proof}

\section{Proof of Theorem~\ref{thm:rates}}

Let $Z=(X,Y)$ on $S=\T\times\T$, with stationary law $\pi\otimes\pi$ (again written $\pi$), and
\[
B_0=D_0=\{(x,y):|x-y|\le\rho\},\qquad D^X_j=A_j\times\T,\qquad D^Y_j=\T\times A_j ,
\]
ordered as $D_0,D^X_1,\dots,D^X_k,D^Y_1,\dots,D^Y_k$ to match the priority convention.

\paragraph{Rates.} $X-Y$ is a rate-$2$ simple random walk on $\T$, started uniformly under $\pi$,
so $\EE_\pi T_{D_0}=\tfrac12\EE_\pi T_{B(0,\rho)}$, i.e.\ the rate of $D_0$ is $2h(\rho)$. The rates
of $D^X_j$ and $D^Y_j$ are $h_j$. By Lemma~\ref{lem:hit}(b) every rate $h'$ in this list satisfies
$th'\le Cg/\log n\to0$, and $\trel/t\le1/(4g)$. Hence, by Lemma~\ref{lem:AB}, each
$a_i=\PP_\pi(T_i\le t)$ equals $t$ times its rate times $1+\eta_i$, where
$\eta:=\max_i|\eta_i|\le\frac1{4g}+Cg/\log n\to0$. In particular
$a=2t\,(h(\rho)+H)(1+O(\eta))$.

\paragraph{Pairs.} We show $\delta=o(a)$. There are finitely many pairs, of three types.

\emph{(i) $D^X_j$ and $D^Y_l$.} $X$ and $Y$ are independent under $\pi$, so the probability that
both are hit by time $t$ is $a^X_ja^Y_l\le Cth_j\cdot a$.

\emph{(ii) $D^X_j$ and $D^X_l$, $j\ne l$} (and likewise for $Y$). By the strong Markov property
at $T^X_j$ and Lemma~\ref{lem:hit}(c), as $\dist(A_j,A_l)\ge n/\log n$,
\[
\PP_\pi(T^X_j\le T^X_l\le t)\le a^X_j\sup_{x\in A_j}\PP_x(T_{A_l}\le t)\le a^X_j\kappa_n ,
\]
so the pair contributes at most $2\kappa_na$.

\emph{(iii) $D_0$ and $D^X_j$} (and likewise for $Y$). Write $D_u=X_u-Y_u$.
Realise the stationary pair as follows: $X_0$ and $D_0$ are independent and uniform; jumps occur
at the points of a rate-$2$ Poisson process; at the $i$-th jump an independent uniform direction
$e_i$ and an independent fair sign $\sigma_i$ are drawn, and $X$ moves by $e_i$ if $\sigma_i=+1$
while $Y$ moves by $e_i$ if $\sigma_i=-1$. Then $D$ jumps by $f_i=\sigma_ie_i$. The directions
$(f_i)$ are i.i.d.\ uniform and independent of $(\sigma_i)$, so the path of $D$ is a function of
$(D_0,\text{jump times},(f_i))$, which is independent of $(X_0,(\sigma_i))$. Since
$X_u=X_0+\sum_{i\le N_u}f_i\one_{\{\sigma_i=+1\}}$ and $X_0$ is uniform and independent of
everything else, $X_{T_0}$ is uniform and independent of the path of $D$. By the strong Markov
property of the pair at $T_0$,
\begin{equation}\label{eq:F}
\PP_\pi\bigl(T_0\le2t,\ X\text{ visits }A_j\text{ during }[T_0,T_0+2t]\bigr)
=\PP_\pi(T_0\le2t)\,\PP_\pi(T_{A_j}\le2t).
\end{equation}
Now let $F=\{\exists\,u\le v\le2t:\ X_u\in A_j,\ D_v\in B(0,\rho)\}$ and
$\tilde F=\{\exists\,v\le u\le2t:\ D_v\in B(0,\rho),\ X_u\in A_j\}$. The stationary pair is
reversible, so its time reversal on $[0,2t]$ has the same law and $\PP_\pi(F)=\PP_\pi(\tilde F)$.
The event $\tilde F$ is contained in the event in \eqref{eq:F}, and
$\{T^X_j\le t,\ T_0\le t\}\subset F\cup\tilde F$. With Lemma~\ref{lem:AB},
\[
\PP_\pi(T^X_j\le t,\ T_0\le t)\le2\,(2t+\trel)^2\,2h(\rho)\,h_j\le Cth(\rho)\cdot a .
\]
Altogether $\delta\le C(\kappa_n+g/\log n)\,a=o(a)$.

\paragraph{Short times.} By Lemma~\ref{lem:AB},
$\PP_\pi(T\le s)\le(s+\trel)\cdot2(h(\rho)+H)\le C(s/t)\,a=Cg^{-1/2}a$.
By Lemma~\ref{lem:mix}, $d(s)\le(1-c^2)^{\lfloor g^{1/2}\rfloor}\to0$.

\paragraph{Conclusion from stationarity.} Lemma~\ref{lem:comp} now gives
$\bar p_i=a_i/a+o(1)$, and $a_i/a$ equals the ratio of the corresponding rates up to a factor
$1+O(\eta)$. Thus $\bar p_0=h(\rho)/(h(\rho)+H)+o(1)$ and
$\bar p^X_j=\bar p^Y_j=h_j/(2(h(\rho)+H))+o(1)$.

\paragraph{Admissible starts.} Let $z_0=(x_0,y_0)$ be admissible. By Lemma~\ref{lem:hit}(c),
$\PP_{x_0}(T_{A_j}\le s)\le\kappa_n$, and the same holds for $y_0$. The difference $X-Y$ is a rate-$2$ walk
started at distance at least $n/\log n$ from $B(0,\rho)$, so
$\PP_{z_0}(T_0\le s)=\PP_{x_0-y_0}(T_{B(0,\rho)}\le2s)\le\kappa_n$ for a rate-$1$ walk, since
$2s\le2t$. Hence $\PP_{z_0}(T\le s)\le(2k+1)\kappa_n$, and Lemma~\ref{lem:start} gives
$|f_i(z_0)-\bar p_i|\le d(s)+(2k+1)\kappa_n\to0$. All bounds depend on the configuration only
through $k$ and $\alpha_*$. \qed

\begin{remark}
In the point-like case, $\EE_\pi T_{\{a\}}$ is the same for every site $a$, so the three rates are
$2h(0),h(0),h(0)$ exactly and the limit $1/(1+k)$ uses only Lemma~\ref{lem:hit}(b),(c).
\end{remark}

\section{Proof of Theorem~\ref{thm:fate}}

\paragraph{Single-walk competition.} Apply Section~2 to one walk with the sets $A_1,\dots,A_k$.
The pair bound (ii), the short-time bound and the mixing bound hold verbatim, so with
$\varphi_l(y)=\PP_y(\text{the first hole touched is }A_l)$ and $\bar\varphi_l=\EE_\pi\varphi_l$,
\begin{equation}\label{eq:single}
\bar\varphi_l=\frac{h_l}{H}+o(1),\qquad
|\varphi_l(y)-\bar\varphi_l|\le d(s)+k\kappa_n\quad\text{if }\dist(y,A_j)\ge\frac n{\log n}\ \forall j .
\end{equation}

\paragraph{Position of the free defect.} Let $E^X_i$ be the event that the first event is the
anchoring of $X$ at $A_i$. Let $\tau^X=\min_jT^X_j$, which is a function of the path of $X$ only.
On $E^X_i$ we have $T=\tau^X$, and $Y_T$ lies in no hole. Let
$\mathcal N=\{y:\dist(y,A_j)<n/\log n\text{ for some }j\}$, so that
$\pi(\mathcal N)\le Ck(n/\log n+n^{\alpha_*})^2/n^2\to0$. Since $\tau^X$ is independent of $Y$,
\[
\PP(E^X_i,\ Y_T\in\mathcal N)\le\PP(Y_{\tau^X}\in\mathcal N)
\le\PP_{x_0}(\tau^X\le s)+\sup_{u\ge s}\PP_{y_0}(Y_u\in\mathcal N)
\le k\kappa_n+\pi(\mathcal N)+d(s)\to0 .
\]

\paragraph{Second phase.} By the strong Markov property at $T$, on $E^X_i$ the free defect $Y$
continues as a walk from $Y_T$, so
$\PP(E^X_i,\ Y\text{ next touches }A_l)=\EE[\one_{E^X_i}\varphi_l(Y_T)]$. By \eqref{eq:single} and
the previous display,
\[
\PP(E^X_i,\ Y\text{ next touches }A_l)=\PP(E^X_i)\,\frac{h_l}{H}+o(1)
=\frac{h_ih_l}{2H\,(h(\rho)+H)}+o(1),
\]
using Theorem~\ref{thm:rates}. The same holds with the roles of $X$ and $Y$ exchanged.

\paragraph{Summation.} $S_{\{i,l\}}$ occurs in four disjoint ways ($X$ or $Y$ anchored first, at
$A_i$ or at $A_l$), so $\PP(S_{\{i,l\}})=\frac{2h_ih_l}{H(h(\rho)+H)}+o(1)$. Termination occurs by
annihilation or by both defects reaching the same hole:
$\PP(\text{termination})=\frac{h(\rho)}{h(\rho)+H}+\frac{\sum_ih_i^2}{H(h(\rho)+H)}+o(1)$.
By Theorem~\ref{thm:asym}, $h_j/h(\rho)\to w_j/w_0$, which gives the stated limits; they sum to
$1$. For $k=1$, $W^2=Q$. For point-like core and holes, $w_0=w_j=1$, $W=Q=k$, and
$(W^2-Q)/(W(1+W))=(k-1)/(k+1)$. \qed

\section{Remarks}

\begin{remark}[Numerical check]
Monte Carlo simulation on tori with $n=32,64,128$ and $k=1,2,4$ point holes (3{,}000--6{,}000
runs each, stationary starts) gives $\PP(\text{annihilation first})$ and
$\PP(\text{sustained pair})$ within $0.03$ of $1/(1+k)$ and $(k-1)/(k+1)$; for $k=1$ no
sustained pair occurs.
\end{remark}

\begin{remark}[Scope]
(i) The error terms are $O((\log n)^{-1/4})$ from the choice of $g$ and are not optimised.
(ii) The proof uses the torus only through Lemmas~\ref{lem:hk} and \ref{lem:mix} and through
translation invariance (for the identity of rates and for step (iii), which uses the group
structure). (iii) Partial anchoring (a hole anchors with probability $q<1$ per visit), creation of
new pairs, detachment, interaction between defects and anisotropic conduction are not covered by
these theorems.
\end{remark}

\end{document}